\documentclass[11pt]{amsart}

\usepackage{color, amsmath,amssymb, amsfonts, amstext,amsthm, latexsym}

   \usepackage{latexsym, oldlfont}

\newcommand{\R}{{\mathbb R}}

\newcommand{\beq}{\begin{equation}}
\newcommand{\eeq}{\end{equation}}
\newcommand{\Del}{\Delta}

\newcommand{\nab}{\nabla}
\newcommand{\lam}{\lambda}

\newcommand{\Ome}{\Omega}
\newcommand{\ome}{\omega}
\newcommand{\patl}{\partial}

\newcommand{\al}{\alpha}

\title[Navier-Stokes Equations]{Symplectic Representation and Turbulent
Global Solutions of Incompressible Navier-Stokes Equations in $\R^3$}

\author{Yongqian Han}

\date{}

\address{
Institute of Applied Physics and Computational Mathematics,
Beijing 100088, China}
\address{
National Key Laboratory of Computational Physics,
Beijing 100088, China}
\email{han\_yongqian@iapcm.ac.cn}



\begin{document}

\begin{abstract} The incompressible Navier-Stokes equations are
considered. We find that there exist infinite non-trivial solutions of
static Euler equations. Moreover there exist random solutions of
static Euler equations. Provided Reynolds number is large enough and
time variable $t$ goes to infinity, these random solutions of static
Euler equations are the path limits of corresponding Navier-Stokes
flows. But the double limit of these Navier-Stokes flows do not exist.
Therefore these solutions are called turbulent solutions.
\end{abstract}

\maketitle

~ ~ ~ ~
%{\bf Chinese Library Classification}: O175.29

\section{Introduction}
\setcounter{equation}{0}

The Navier-Stokes equations in $\R^3$ are given by
\beq\label{NS1}
u_t-\nu\Del u+(u\cdot\nab)u+\nab P=0,
\eeq
\beq\label{NS2}
\nab\cdot u=0,
\eeq
where $u=u(t,x)=\big(u^1(t,x),u^2(t,x),u^3(t,x)\big)$ and $P=P(t,x)$
stand for the unknown velocity vector field of fluid and its
pressure, $\nu>0$ is the coefficient of viscosity. Here $x=(x_1,x_2,
x_3)$, $\nab=(\patl_1,\patl_2,\patl_3)$ and $\patl_j=\patl_{x_j}$
($j=1,2,3$).

Taking curl with equation \eqref{NS1}, we have
\beq\label{NS-G-ome}
\ome_t-\nu\Del\ome+(u\cdot\nab)\ome-(\ome\cdot\nab)u=0.
\eeq
Here vorticity $\ome=(\ome^1,\ome^2,\ome^3)$ and
\beq\label{Def-G-ome}\begin{split}
\ome(t,x)=&\nab\times u(t,x).\\
\end{split}\eeq

The velocity vector $u$ is solution of equations \eqref{NS1}
\eqref{NS2} provided $(u,\ome)$ satisfies equation \eqref{NS-G-ome}.

The Fourier transformation of $u(t,x)$ with respect to $x$ denotes by
$\hat{u}(t,\xi)$. Then equation \eqref{NS2} is equivalent to
\beq\label{NS2-F}
\xi_j\hat{u}^j=\xi_1\hat{u}^1+\xi_2\hat{u}^2+\xi_3\hat{u}^3=0.
\eeq
It is that $\hat{u}(t,\xi)$ is perpendicular to $\xi$. Denote by $\xi
\bot\hat{u}$. Equivalently $\hat{u}(t,\xi)\in T_{\xi}{\mathbb{S}}^2$
for any $\xi\in{\mathbb{S}}^2$.

For any vector $A\in\R^3-\{0\}$, let
\beq\label{Sym-Rep-12}
u(t,x)=\{A\times\nab\}\phi(t,x)+\{(A\times\nab)\times\nab\}\psi(t,x),
\eeq
where $\phi$ and $\psi$ are real scalar functions. Then $u$ satisfies
the equation \eqref{NS2}. We call that the formulation
\eqref{Sym-Rep-12} is (1,2)-symplectic representation of
velocity vector $u$.

Now
\beq\label{Def-SR12-ome}\begin{split}
\ome(t,x)=&\nab\times u(t,x)\\
=&-\{(A\times\nab)\times\nab\}\phi(t,x)+(A\times\nab)\Del\psi(t,x).\\
\end{split}\eeq
Thanks the following observations
\beq\label{SR12-phi}\begin{split}
(A\times\nab)\cdot u(t,x)&=(A\times\nab)\cdot(A\times\nab)\phi(t,x)\\
&=\big\{(A\cdot A)\Del-(A\cdot\nab)^2\big\}\phi(t,x),\\
\end{split}\eeq
\beq\label{SR12-psi}\begin{split}
(A\times\nab)\cdot \ome(t,x)&=(A\times\nab)\cdot(A\times\nab)\Del\psi(t,x)\\
&=\big\{(A\cdot A)\Del-(A\cdot\nab)^2\big\}\Del\psi(t,x),\\
\end{split}\eeq
taking scalar product of equation \eqref{NS1} with $A\times\nab$,
we have
\beq\label{NS-SR12-phi}
\big\{(A\cdot A)\Del-(A\cdot\nab)^2\big\}\{\phi_t-\nu\Del\phi\}
+(A\times\nab)\cdot\{(u\cdot\nab)u\}=0.
\eeq
And taking scalar product of \eqref{NS-G-ome} with $A\times\nab$,
we derive
\beq\label{NS-SR12-psi}
\big\{(A\cdot A)\Del-(A\cdot\nab)^2\big\}\Del\{\psi_t-\nu\Del\psi\}
+(A\times\nab)\cdot\{(u\cdot\nab)\ome-(\ome\cdot\nab)u\}=0.
\eeq

Therefore $(\phi,\psi)$ satisfies equations \eqref{NS-SR12-phi}
\eqref{NS-SR12-psi} provided $u$ satisfies equations \eqref{NS1}
\eqref{Sym-Rep-12}.

There is a large literature (see \cite{Con01,Te01} and therein
references) studying the incompressible Navier-Stokes equations. Here
we construct some turbulent global solutions of equations \eqref{NS1}
\eqref{NS2} by using \eqref{Sym-Rep-12} \eqref{NS-SR12-phi}
\eqref{NS-SR12-psi}, provided $(\phi,\psi)$ is radial symmetric and
cylindrical symmetric respectively.

%%%%%%%%%%%%%%%%%%%
\section{Radial Symplectic Representation}
\setcounter{equation}{0}
%%%%%%%%%%%%%%%%%%%%%%%%%%%%%%%%%%%

In this section, we assume that $\phi$ and $\psi$ are radial symmetric
functions with respect to space variable $x\in\R^3$. It is that
$\phi(t,x)=\phi(t,r)$, $\psi(t,x)=\psi(t,r)$  and $r^2=x_1^2+x_2^2
+x_3^2$. Inserting $\big(\phi(t,r),\psi(t,r)\big)$ into
\eqref{NS-SR12-phi} \eqref{NS-SR12-psi}, by long long straightforward
calculation, we find that $u$ is the solution of Navier-Stokes
equations \eqref{NS1} \eqref{NS2}, provided
\beq\label{Sol-RSR12-NS}\begin{split}
u(t,x)=&e^{-\nu\lam^2t}U(x)\\
U(x)=&\{A\times\nab\}\Phi(r)+\{(A\times\nab)\times\nab\}\Psi(r),
\end{split}\eeq
\beq\label{RSR12-phL}\begin{split}
\Phi=&\lam\Psi,\\
\end{split}\eeq
\beq\label{RSR12-psL}\begin{split}
-\Del\Psi=&\lam^2\Psi.\\
\end{split}\eeq
Here vector $A\in\R^3-\{0\}$.

Moreover $U(x)$ is the solution of static Euler equations
\beq\label{STEE-RS}\begin{split}
&(U\cdot\nab)U+\nab \tilde{P}=0,\\
&\nab\cdot U=0,\\
\end{split}\eeq
provided $\big(\Phi(r),\Psi(r)\big)$ satisfies \eqref{RSR12-phL}
\eqref{RSR12-psL}.

It is well-known that the equation \eqref{STEE-RS} is equivalent to
\beq\label{STEE-curl}\begin{split}
&(U\cdot\nab)\Ome-(\Ome\cdot\nab)U
=\sum_{j=1}^3\patl_j(U^j\Ome-\Ome^jU)=0,\\
\end{split}\eeq
where vorticity $\Ome=\nab\times U$.

For simplicity let us take $A=e_3=(0,0,1)$, then
\beq\label{RSR12-u-ome}\begin{split}
U(x)=&\{e_3\times\nab\}\Phi(r)+\{(e_3\times\nab)\times\nab\}\Psi(r)\\
=&(-\patl_2\Phi,\patl_1\Phi,0)+\big(\patl_1\patl_3\Psi,
\patl_2\patl_3\Psi,-(\patl_1^2+\patl_2^2)\Psi\big),\\
\Ome(x)=&\nab\times U(x)\\
=&\big(-\patl_1\patl_3\Phi,-\patl_2\patl_3\Phi,
(\patl_1^2+\patl_2^2)\Phi\big)+(-\patl_2\Del\Psi,\patl_1\Del\Psi,0).
\end{split}\eeq
It is easy to verify that equation \eqref{STEE-curl} is satisfied
provided $\big(\Phi(r),\Psi(r)\big)$ satisfies \eqref{RSR12-phL}
\eqref{RSR12-psL}. Indeed we have
\beq\label{STEE-curl-et}\begin{split}
&U^j\Ome^k-\Ome^jU^k=0,\;\;\;\;\;\;j,k=1,2,3.\\
\end{split}\eeq
By straightforward calculation, we have
\beq\label{uo-21}\begin{split}
&\{U^2\Ome^1-\Ome^2U^1\}\Big|_{\Phi=\lam\Psi}\\
=&\{(\patl_1\Phi+\patl_2\patl_3\Psi)(-\patl_1\patl_3\Phi-\patl_2\Del\Psi)\\
&-(-\patl_2\patl_3\Phi+\patl_1\Del\Psi)(-\patl_2\Phi+\patl_1\patl_3\Psi)\}\Big|_{\Phi=\lam\Psi}\\
=&-(\lam^2\patl_1\Psi+\patl_1\Del\Psi)\patl_1\patl_3\Psi
-(\lam^2\patl_2\Psi+\patl_2\Del\Psi)\patl_2\patl_3\Psi\\
&-\lam(\patl_1\Psi\patl_2\Del\Psi-\patl_2\Psi\patl_1\Del\Psi),\\
\end{split}\eeq

\beq\label{uo-31}\begin{split}
&\{U^3\Ome^1-\Ome^3U^1\}\Big|_{\Phi=\lam\Psi}\\
=&\{-(\patl_1^2+\patl_2^2)\Psi(-\patl_1\patl_3\Phi-\patl_2\Del\Psi)
-(\patl_1^2+\patl_2^2)\Phi(-\patl_2\Phi+\patl_1\patl_3\Psi)\}\Big|_{\Phi=\lam\Psi}\\
=&(\patl_1^2+\patl_2^2)\Psi(\patl_2\Del\Psi+\lam^2\patl_2\Psi),\\
\end{split}\eeq

\beq\label{uo-12}\begin{split}
&\{U^1\Ome^2-\Ome^1U^2\}\Big|_{\Phi=\lam\Psi}\\
=&\{(-\patl_2\Phi+\patl_1\patl_3\Psi)(-\patl_2\patl_3\Phi+\patl_1\Del\Psi)\\
&-(-\patl_1\patl_3\Phi-\patl_2\Del\Psi)(\patl_1\Phi+\patl_2\patl_3\Psi)\}\Big|_{\Phi=\lam\Psi}\\
=&(\lam^2\patl_2\Psi+\patl_2\Del\Psi)\patl_2\patl_3\Psi
+(\lam^2\patl_1\Psi+\patl_1\Del\Psi)\patl_1\patl_3\Psi\\
&+\lam(\patl_1\Psi\patl_2\Del\Psi-\patl_2\Psi\patl_1\Del\Psi),\\
\end{split}\eeq

\beq\label{uo-32}\begin{split}
&\{U^3\Ome^2-\Ome^3U^2\}\Big|_{\Phi=\lam\Psi}\\
=&\{-(\patl_1^2+\patl_2^2)\Psi(-\patl_2\patl_3\Phi+\patl_1\Del\Psi)
-(\patl_1^2+\patl_2^2)\Phi(\patl_1\Phi+\patl_2\patl_3\Psi)\}\Big|_{\Phi=\lam\Psi}\\
=&-(\patl_1^2+\patl_2^2)\Psi(\lam^2\patl_1\Psi+\patl_1\Del\Psi),\\
\end{split}\eeq

\beq\label{uo-13}\begin{split}
&\{U^1\Ome^3-\Ome^1U^3\}\Big|_{\Phi=\lam\Psi}\\
=&\{(-\patl_2\Phi+\patl_1\patl_3\Psi)(\patl_1^2+\patl_2^2)\Phi
+(-\patl_1\patl_3\Phi-\patl_2\Del\Psi)(\patl_1^2+\patl_2^2)\Psi\}\Big|_{\Phi=\lam\Psi}\\
=&-(\lam^2\patl_2\Psi+\patl_2\Del\Psi)(\patl_1^2+\patl_2^2)\Psi,\\
\end{split}\eeq

\beq\label{uo-23}\begin{split}
&\{U^2\Ome^3-\Ome^2U^3\}\Big|_{\Phi=\lam\Psi}\\
=&\{(\patl_1\Phi+\patl_2\patl_3\Psi)(\patl_1^2+\patl_2^2)\Phi
+(-\patl_2\patl_3\Phi+\patl_1\Del\Psi)(\patl_1^2+\patl_2^2)\Psi\}\Big|_{\Phi=\lam\Psi}\\
=&(\lam^2\patl_1\Psi+\patl_1\Del\Psi)(\patl_1^2+\patl_2^2)\Psi.\\
\end{split}\eeq

It is obvious that equation \eqref{STEE-curl-et} is satisfied provided
radial function $\Psi(r)$ satisfies equation \eqref{RSR12-psL}.

Here $u(t,x)=e^{-\nu\lam^2t}U(x)$ and $\ome(t,x)=\nab\times u(t,x)=
e^{-\nu\lam^2t}\Ome(x)$. It is easy to verify that $(u,\ome)$
satisfies \eqref{NS-G-ome}.

The equation \eqref{RSR12-psL} is equivalent to
\beq\label{RSR12-psL1}\begin{split}
-\patl_r^2(r\Psi)=&\lam^2r\Psi.\\
\end{split}\eeq
There exist infinite solutions
\beq\label{Sol-RSR-ps}\begin{split}
&\Psi(r)=\Psi_{\lam\al\beta}(r)=\al\frac1r\sin(\lam r)
+\beta\frac1r\cos(\lam r)\\
\end{split}\eeq
of equation \eqref{RSR12-psL1} for any real constants $\lam,\al$ and
$\beta$.

Therefore there exist solution class ${\mathbb X}_E$ of static Euler
equations and solution class ${\mathbb X}_{NS}$ of Navier-Stokes
equations. Here
\beq\label{SC-RS-E}\begin{split}
&{\mathbb X}_E=\cup_{\lam\in\R}X_{e\lam},\\
&X_{e\lam}=\mbox{span}\big\{\lam(e_3\times\nab)\big(\frac1r\sin(\lam r)
\big)+\{(e_3\times\nab)\times\nab\}\big(\frac1r\sin(\lam r)\big),\\
&\hspace{23mm}\lam(e_3\times\nab)\big(\frac1r\cos(\lam r)\big)
+\{(e_3\times\nab)\times\nab\}\big(\frac1r\cos(\lam r)\big)\big\},\\
\end{split}\eeq
\beq\label{SC-RS-NS}\begin{split}
&{\mathbb X}_{NS}=\cup_{\lam\in\R}X_{ns\lam},\\
&X_{ns\lam}=\mbox{span}\big\{u_{\lam}^s,u_{\lam}^c\big\}\\
&u_{\lam}^s=\lam e^{-\nu\lam^2t}(e_3\times\nab)\big(
\frac1r\sin(\lam r)\big)+e^{-\nu\lam^2t}\{(e_3\times\nab)\times\nab\}
\big(\frac1r\sin(\lam r)\big),\\
&u_{\lam}^c=\lam e^{-\nu\lam^2t}(e_3\times\nab)\big(\frac1r\cos(\lam r)
\big)+e^{-\nu\lam^2t}\{(e_3\times\nab)\times\nab\}
\big(\frac1r\cos(\lam r)\big),\\
\end{split}\eeq
where $X_{e\lam}$ and $X_{ns\lam}$ are linear spaces for any
$\lam\in\R$.

There also exist non-trivial solutions of equation \eqref{RSR12-psL}
with boundary condition.

Let $R_0>0$ and ball $B_0=\{r<R_0\}\subset\R^3$. We consider Dirichlet
eigenvalue problem
\beq\label{L-eq3}\begin{split}
-(\patl_1^2+\patl_2^2+\patl_3^2)\Psi=&\lam^2\Psi,\\
\end{split}\eeq
\beq\label{DBC03}
\Psi|_{r=R_0}=0.
\eeq
The eigenvalue $\lam$ of problem \eqref{L-eq3} \eqref{DBC03} is
so-called Dirichlet eigenvalue. By spectrum theory of compact
operator $(-\Del)^{-1}$, there exist countable infinite Dirichlet
eigenvalues
\[
0<\lam_1\le\cdots\le\lam_n\le\cdots
\]
Each Dirichlet eigenvalue is real, and may be geometrically simple or
finite.

Let $\rho$ be orthogonal transformation in $\R^3$, and $\Psi_j(x)$ be
the eigenfunction corresponding to Dirichlet eigenvalue $\lam_j$. Then
$\Psi_j(\rho x)$ is also eigenfunction, and $\Psi_j(x)$ is radial
function provided $\lam_j$ is simple. Therefore $\lam_1=\frac{\pi}
{R_0}$ and
\beq\label{DEF-1}\begin{split}
&\Psi_1(x)=\Psi_1(r)=\al\frac1r\sin(\lam_1 r),
\end{split}\eeq
since the first eigenvalue $\lam_1$ is simple and $\Psi_1(x)\not=0$
for any $x\in B_0$.

It is obvious that $\Psi_1\in H^1_0(B_0)\cap H^2(B_0)$. By the regular
theory of elliptic equations, we have $\Psi_1\in H^{\infty}_{loc}
(\R^3)$. Since the regularity of $\Psi_1$ is independent of $\lam_1$,
we derive
\[
\Psi_{\lam\al 0}(r)=\Psi_{\lam\al\beta}{\big|}_{\beta=0}(r)
=\al\frac1r\sin(\lam r)\in H^{\infty}_{loc}(\R^3).
\]

Therefore there exist solution class ${\mathbb X}^s_E$ of static Euler
equations and solution class ${\mathbb X}^s_{NS}$ of Navier-Stokes
equations. Here
\beq\label{SC-RS-ES}\begin{split}
&{\mathbb X}^s_E=\cup_{\lam\in\R}X^s_{e\lam}\subset
H^{\infty}_{loc}(\R^3),\\
&X^s_{e\lam}=\mbox{span}\big\{\lam(e_3\times\nab)\big(\frac1r\sin(\lam r)
\big)+\{(e_3\times\nab)\times\nab\}\big(\frac1r\sin(\lam r)\big)\big\},\\
\end{split}\eeq
\beq\label{SC-RS-NSS}\begin{split}
&{\mathbb X}^s_{NS}=\cup_{\lam\in\R}X^s_{ns\lam}\subset
H^{\infty}_{loc}(\R^3),\\
&X^s_{ns\lam}=\mbox{span}\big\{\lam e^{-\nu\lam^2t}(e_3\times\nab)\big(
\frac1r\sin(\lam r)\big)\\
&\hspace{23mm}+e^{-\nu\lam^2t}\{(e_3\times\nab)\times\nab\}
\big(\frac1r\sin(\lam r)\big)\big\},\\
\end{split}\eeq
where $X^s_{e\lam}\subset H^{\infty}_{loc}(\R^3)$ and $X^s_{ns\lam}
\subset H^{\infty}_{loc}(\R^3)$ are linear spaces for any $\lam\in\R$.

%%%%%%%%%%%%%%%%%%%
\section{Cylindrical Symplectic Representation}
\setcounter{equation}{0}
%%%%%%%%%%%%%%%%%%%%%%%%%%%%%%%%%%%

In this section, we assume that $A=e_3=(0,0,1)$, $\phi$ and $\psi$
are cylindrical symmetric functions with respect to space variable
$x\in\R^3$. It is that $\phi(t,x)=\phi(t,r,x_3)$, $\psi(t,x)=\psi
(t,r,x_3)$  and $r^2=x_1^2+x_2^2$.

Assume that symplectic representation $(\phi,\psi)$ is independent of
$x_3$. Inserting $\big(\phi(t,r),$ $\psi(t,r)\big)$ into
\eqref{NS-SR12-phi} \eqref{NS-SR12-psi}, by straightforward
calculation, we find that there exist solution class
${\mathbb X}^{2R}_E$ of static Euler equations and
solution class ${\mathbb X}^{2R}_{NS}$ of Navier-Stokes equations.
Here
\beq\label{SC-CSR12-E2r}\begin{split}
{\mathbb X}^{2R}_E=\big\{&(A\times\nab)\Phi(r)
+\{(A\times\nab)\times\nab\}\Psi(r)\\
&\big|\mbox{$\Phi(r)$ and $\Psi(r)$ are any functions}\big\},
\end{split}\eeq
\beq\label{SC-CSR12-NS2r}\begin{split}
{\mathbb X}^{2R}_{NS}=\big\{&(A\times\nab)e^{-t\nu\Del}\Phi(r)
+\{(A\times\nab)\times\nab\}e^{-t\nu\Del}\Psi(r)\\
&\big|\mbox{$\Phi(r)$ and $\Psi(r)$ are any functions}\big\}.
\end{split}\eeq
Indeed
\[
U=U(x_1,x_2)=(A\times\nab)\Phi(r)
+\{(A\times\nab)\times\nab\}\Psi(r)
\]
is the solution of static Euler equations \eqref{STEE-RS}, and
\[
u=u(t,x_1,x_2)=
(A\times\nab)e^{-t\nu\Del}\Phi(r)+\{(A\times\nab)\times\nab\}
e^{-t\nu\Del}\Psi(r)
\]
is the solution of Navier-Stokes equations \eqref{NS1} \eqref{NS2}.

Assume that symplectic representation $(\phi,\psi)$ depends on $x_3$.
Inserting $\big(\phi(t,r,x_3),$ $\psi(t,r,x_3)\big)$ into
\eqref{NS-SR12-phi} \eqref{NS-SR12-psi}, by long long straightforward
calculation, we find that $u$ is the solution of Navier-Stokes
equations \eqref{NS1} \eqref{NS2}, provided
\beq\label{Sol-CSR12-NS}\begin{split}
u(t,x)=&e^{-\nu\lam^2t}U(x)\\
U(x)=&\{A\times\nab\}\Phi(r,x_3)+\{(A\times\nab)\times\nab\}\Psi(r,x_3),
\end{split}\eeq
\beq\label{CSR12-phL}\begin{split}
\Phi=&\lam\Psi,\\
\end{split}\eeq
\beq\label{CSR12-psL}\begin{split}
-\Del\Psi=&\lam^2\Psi.\\
\end{split}\eeq

Moreover
\beq\label{Sol-CSR12-EE}\begin{split}
U(x)=&\{e_3\times\nab\}\Phi(r,x_3)
+\{(e_3\times\nab)\times\nab\}\Psi(r,x_3),\\
\Ome(x)=&\nab\times U(x),
\end{split}\eeq
and $U(x)$ is the solution of static Euler equations \eqref{STEE-RS}
provided $\big(\Phi(r,x_3),\Psi(r,x_3)\big)$ satisfies
\eqref{CSR12-phL} \eqref{CSR12-psL}. Indeed the equations
\eqref{RSR12-u-ome}--\eqref{uo-23} are also satisfied by $(U,\Ome)$
defined in \eqref{Sol-CSR12-EE}.

Here $u(t,x)=e^{-\nu\lam^2t}U(x)$ and $\ome(t,x)=\nab\times u(t,x)=
e^{-\nu\lam^2t}\Ome(x)$. It is easy to verify that $(u,\ome)$
satisfies \eqref{NS-G-ome}.

There exist non-trivial solutions of equation \eqref{CSR12-psL}
with boundary conditions. Now we construct some solutions of
\eqref{CSR12-psL}.

Let $R_0>0$ and disc $D_0=\{r<R_0\}\subset\R^2$. We consider Dirichlet
eigenvalue problem
\beq\label{L-eq2}\begin{split}
-(\patl_1^2+\patl_2^2)W=&\kappa^2W,\\
\end{split}\eeq
\beq\label{DBC02}
W|_{r=R_0}=0.
\eeq
By spectrum theory of compact operator $(-\Del)^{-1}$, there exist
countable infinite Dirichlet eigenvalues
\[
0<\kappa_1\le\cdots\le\kappa_n\le\cdots
\]
Each Dirichlet eigenvalue is real, and may be geometrically simple or
finite.

Let $\tilde{\rho}$ be orthogonal transformation in $\R^2$, and $W_j(x)$ be
the eigenfunction corresponding to Dirichlet eigenvalue $\kappa_j$. Then
$W_j(\tilde{\rho}x)$ is also eigenfunction, and $W_j(x)$ is radial
function provided $\kappa_j$ is simple. Therefore $W_1(x)$ is radial
since the first eigenvalue $\kappa_1$ is simple.

It is obvious that $W_j\in H^1_0(D_0)\cap H^2(D_0)$. By the regular
theory of elliptic equations, we have $W_j\in H^{\infty}(D_0)$ for
any $j=1,2,\cdots$.

Therefore there exist solutions
\beq\label{Sol-CSR-psD}\begin{split}
&\lam_j^2=\kappa_j^2+\eta^2,\;\;j=1,2,\cdots\\
&\Psi_j(r,x_3)=\al W_j(r)\sin(\eta x_3)+\beta W_j(r)\cos(\eta x_3)\\
\end{split}\eeq
of problem \eqref{CSR12-psL} provided $W_j(x)$ is radial. Here
$\al$, $\beta$ and $\eta$ are any real constants.

Thus there exist solution class ${\mathbb Y}_E$ of static Euler
equations and solution class ${\mathbb Y}_{NS}$ of Navier-Stokes
equations. Here
\beq\label{SC-RS-Ey}\begin{split}
&{\mathbb Y}_E=\cup^{\infty}_{j=1}\cup_{\eta\in\R}Y^e_{j\eta},\\
&Y^e_{j\eta}=\mbox{span}\big\{\lam_j(A\times\nab)\big(W_j(r)\sin(\eta x_3)
\big)+\{(A\times\nab)\times\nab\}\big(W_j(r)\sin(\eta x_3)\big),\\
&\hspace{22mm}\lam_j(A\times\nab)\big(W_j(r)\cos(\eta x_3)\big)
+\{(A\times\nab)\times\nab\}\big(W_j(r)\cos(\eta x_3)\big)\big\},\\
\end{split}\eeq
\beq\label{SC-CS-NSy}\begin{split}
{\mathbb Y}_{NS}=&\cup^{\infty}_{j=1}\cup_{\eta\in\R}Y^{ns}_{j\eta},\\
Y^{ns}_{j\eta}=&\mbox{span}\big\{u^s_{j\eta},u^c_{j\eta}\big\},\\
u^s_{j\eta}=&\lam_je^{-\nu\lam^2_jt}(A\times\nab)
\big(W_j(r)\sin(\eta x_3)\big)\\
&+e^{-\nu\lam^2_jt}\{(A\times\nab)\times\nab\}
\big(W_j(r)\sin(\eta x_3)\big),\\
u^c_{j\eta}=&\lam_je^{-\nu\lam^2_jt}(A\times\nab)
\big(W_j(r)\cos(\eta x_3)\big)\\
&+e^{-\nu\lam^2_jt}\{(A\times\nab)\times\nab\}
\big(W_j(r)\cos(\eta x_3)\big),\\
\end{split}\eeq
provided $W_j(x)$ is radial. It is obvious that $Y^e_{j\eta}$ and
$Y^{ns}_{j\eta}$ are linear spaces for any $j=1,2,\cdots$ and
$\eta\in\R$.

We observe that
\[
\Del_rV=(\patl_r^2+\frac1r\patl_r)V=\frac1r\patl_r(r\patl_r V).
\]
Let $0<R_1<R_2$ and ring ${\cal R}=\{R_1<r<R_2\}\subset\R^2$. We
consider equation
\beq\label{SL-eq}\begin{split}
-\patl_r(r\patl_r V)=\zeta^2rV\\
\end{split}\eeq
with separated or coupled boundary condition.

The separated boundary condition is defined by
\beq\label{SL-BC1}\begin{split}
& m_{11}V(R_1)+m_{12}\{r\patl_rV(r)\}\big|_{r=R_1}=0,\\
& m_{21}V(R_2)+m_{22}\{r\patl_rV(r)\}\big|_{r=R_2}=0,\\
\end{split}\eeq
where $m_{ij}\in\R$, $(m_{11},m_{12})\not=(0,0)$ and
$(m_{21},m_{22})\not=(0,0)$.

The coupled boundary condition is defined by
\beq\label{SL-BC2}\begin{split}
& Y(r)=\left(\begin{array}{c} V(r)\\r\patl_rV(r)\end{array}\right),\\
& Y(R_1)=KY(R_2),\\
\end{split}\eeq
where $2\times2$ matrix $K\in SL_2(\R)$ and $\det K$=1.

By Sturm-Liouville theory (see \cite{Ze05}, Page 72, Theorem 4.3.1),
there exist countable infinite eigenvalues
\[
0<\zeta_1\le\cdots\le\zeta_n\le\cdots
\]
and eigenfunction $V_j(r)$ corresponding to eigenvalue $\zeta_j$ of
equation \eqref{SL-eq} with separated boundary condition
\eqref{SL-BC1} or coupled boundary condition \eqref{SL-BC2}.
Each eigenvalue is real, and may be geometrically simple or double.

For the eigenfunctions $V_j$ of problem \eqref{SL-eq} \eqref{SL-BC1},
by the regular theory of elliptic equations, we have $V_j\in H^1_0
({\cal R})\cap H^{\infty}({\cal R})$ ($j=1,2,\cdots$) provided
$(m_{11},m_{12})=(1,0)$ and $(m_{21},m_{22})=(1,0)$.

Therefore there exist solutions
\beq\label{Sol-CSR-ps-LT}\begin{split}
&\lam_j^2=\zeta_j^2+\eta^2,\;\;j=1,2,\cdots\\
&\Psi_j(r,x_3)=\al V_j(r)\sin(\eta x_3)+\beta V_j(r)\cos(\eta x_3)\\
\end{split}\eeq
of problem \eqref{CSR12-psL}. Here $\al$, $\beta$ and $\eta$
are any real constants.

Thus there exist solution class ${\mathbb Z}_E$ of static Euler
equations and solution class ${\mathbb Z}_{NS}$ of Navier-Stokes
equations. Here
\beq\label{SC-CS-Ez}\begin{split}
&{\mathbb Z}_E=\cup^{\infty}_{j=1}\cup_{\eta\in\R}Z^e_{j\eta},\\
&Z^e_{j\eta}=\mbox{span}\big\{\lam_j(A\times\nab)\big(V_j(r)\sin(\eta x_3)
\big)+\{(A\times\nab)\times\nab\}\big(V_j(r)\sin(\eta x_3)\big),\\
&\hspace{22mm}\lam_j(A\times\nab)\big(V_j(r)\cos(\eta x_3)\big)
+\{(A\times\nab)\times\nab\}\big(V_j(r)\cos(\eta x_3)\big)\big\},\\
\end{split}\eeq
\beq\label{SC-CS-NSz}\begin{split}
{\mathbb Z}_{NS}=&\cup^{\infty}_{j=1}\cup_{\eta\in\R}Z^{ns}_{j\eta},\\
Z^{ns}_{j\eta}=&\mbox{span}\big\{u^s_{j\eta},u^c_{j\eta}\big\},\\
u^s_{j\eta}=&\lam_je^{-\nu\lam^2_jt}(A\times\nab)
\big(V_j(r)\sin(\eta x_3)\big)\\
&+e^{-\nu\lam^2_jt}\{(A\times\nab)\times\nab\}
\big(V_j(r)\sin(\eta x_3)\big),\\
u^c_{j\eta}=&\lam_je^{-\nu\lam^2_jt}(A\times\nab)
\big(V_j(r)\cos(\eta x_3)\big)\\
&+e^{-\nu\lam^2_jt}\{(A\times\nab)\times\nab\}
\big(V_j(r)\cos(\eta x_3)\big),\\
\end{split}\eeq
where $Z^e_{j\eta}$ and $Z^{ns}_{j\eta}$ are linear spaces for any
$j=1,2,\cdots$ and $\eta\in\R$.

%%%%%%%%%%%%%%%%%%%
\section{Conclusion}
\setcounter{equation}{0}
%%%%%%%%%%%%%%%%%%%%%%%%%%%%%%%%%%%

Some interesting and significant phenomena are founded by observing
the solutions in Section 2 and Section 3.

There exist some special solution classes of Navier-Stokes equations,
for instance ${\mathbb X}_{NS}^{2R}$, ${\mathbb X}_{NS}$,
${\mathbb X}^s_{NS}$, ${\mathbb Y}_{NS}$ and ${\mathbb Z}_{NS}$, which
will be called turbulent solution classes.

Let $\mu$ be random number and $\Gamma_{\mu}$ be a path. In the path
$\Gamma_{\mu}$, we have
\beq\label{Ra-Pa}\begin{split}
t\nu=\mu\\
\end{split}\eeq
for any $(\nu,t)\in\Gamma_{\mu}$.

Assume that $u_{ns}(t,x^{\prime})\in{\mathbb X}_{NS}^{2R}$ is the
solution of Navier-Stokes equation, $x^{\prime}=(x_1,x_2)\in\R^2$ and
$r^{\prime}=|x^{\prime}|$. Then there exist $\Phi(x^{\prime})=\Phi(
r^{\prime})$ and $\Psi(x^{\prime})=\Psi(r^{\prime})$ such that
\[\begin{split}
u_{ns}(t,x^{\prime})=&(e_3\times\nab)\Big\{\frac1{4\pi t\nu}
\int_{\R^2}\exp\{-\frac{|y^{\prime}|^2}{4t\nu}\}\Phi(|x^{\prime}
-y^{\prime}|)dy^{\prime}\Big\}\\
&+\{(e_3\times\nab)\times\nab\}\Big\{\frac1{4\pi t\nu}
\int_{\R^2}\exp\{-\frac{|y^{\prime}|^2}{4t\nu}\}\Psi(|x^{\prime}
-y^{\prime}|)dy^{\prime}\Big\}.\\
\end{split}\]
Letting $(\nu,t)\in\Gamma_{\mu}$ and $(\nu,1/t)\rightarrow(0,0)$, we
have path limit
\beq\label{L-XNS2R}\begin{split}
u_e(x^{\prime})=&\lim_{(\nu,t)\in\Gamma_{\mu},(\nu,1/t)\rightarrow
(0,0)}u_{ns}(t,x^{\prime})\\
=&(e_3\times\nab)\Phi_l(r^{\prime})+\{(e_3\times\nab)\times\nab\}
\Psi_l(r^{\prime}),\\
\Phi_l(r^{\prime})=&\Phi_l(x^{\prime})=\frac1{4\pi\mu}\int_{\R^2}\exp\{
-\frac{|y^{\prime}|^2}{4\mu}\}\Phi(|x^{\prime}-y^{\prime}|)dy^{\prime},\\
\Psi_l(r^{\prime})=&\Psi_l(x^{\prime})=\frac1{4\pi\mu}\int_{\R^2}\exp\{
-\frac{|y^{\prime}|^2}{4\mu}\}\Psi(|x^{\prime}-y^{\prime}|)dy^{\prime}.\\
\end{split}\eeq
It is obvious that random $u_e$ is the solution of static Euler
equations and $u_e\in{\mathbb X}_E^{2R}$. Moreover double limit
$\lim_{(\nu,1/t)\rightarrow(0,0)}u_{ns}(t,x^{\prime})$ does not exist.

Assume that $u_{ns}(t,x)\in{\mathbb X}_{NS}$ is the solution of
Navier-Stokes equation. Then there exist linear space $X_{ns\lam}$ and
real constants $\al$, $\beta$ such that $u_{ns}(t,x)\in X_{ns\lam}$ and
\[\begin{split}
u_{ns}(t,x)=&\al\big\{\lam e^{-\nu\lam^2t}(e_3\times\nab)\big(
\frac1r\sin(\lam r)\big)+e^{-\nu\lam^2t}\{(e_3\times\nab)\times\nab\}
\big(\frac1r\sin(\lam r)\big)\big\}\\
&+\beta\big\{\lam e^{-\nu\lam^2t}(e_3\times\nab)\big(\frac1r\cos(\lam r)
\big)+e^{-\nu\lam^2t}\{(e_3\times\nab)\times\nab\}
\big(\frac1r\cos(\lam r)\big)\big\},\\
\end{split}\]
where $x\in\R^3$ and $r=|x|$. Letting $(\nu,t)\in\Gamma_{\mu}$ and
$(\nu,1/t)\rightarrow(0,0)$, we have path limit
\beq\label{L-XNS-sc}\begin{split}
u_e(x)=&\lim_{(\nu,t)\in\Gamma_{\mu},(\nu,1/t)\rightarrow(0,0)}
u_{ns}(t,x)\\
=&\al\big\{\lam e^{-\mu\lam^2}(e_3\times\nab)\big(
\frac1r\sin(\lam r)\big)+e^{-\mu\lam^2}\{(e_3\times\nab)\times\nab\}
\big(\frac1r\sin(\lam r)\big)\big\}\\
&+\beta\big\{\lam e^{-\mu\lam^2}(e_3\times\nab)\big(\frac1r\cos(\lam r)
\big)+e^{-\mu\lam^2}\{(e_3\times\nab)\times\nab\}
\big(\frac1r\cos(\lam r)\big)\big\}.\\
\end{split}\eeq
It is obvious that random $u_e$ is the solution of static Euler
equations and $u_e\in X_{e\lam}$ since $X_{e\lam}$ is linear space.
Moreover double limit $\lim_{(\nu,1/t)\rightarrow(0,0)}u_{ns}(t,x)$
does not exist.

Similarly assume that $u_{ns}(t,x)\in{\mathbb X}^s_{NS}$ (
${\mathbb Y}_{NS}$ or ${\mathbb Z}_{NS}$ respectively) is the solution
of Navier-Stokes equation. Then there exists linear space $X^s_{ns\lam}$
($Y^{ns}_{j\lam}$ or $Z^{ns}_{j\lam}$ resp.) such that $u_{ns}(t,x)\in
X^s_{ns\lam}$ ($Y^{ns}_{j\lam}$ or $Z^{ns}_{j\lam}$ resp.). Since
$X^s_{e\lam}$ ($Y^e_{j\lam}$ or $Z^e_{j\lam}$ resp.) is linear space, we
can prove that there exists the path limit
\[
u_e(x)=\lim_{(\nu,t)\in\Gamma_{\mu},(\nu,1/t)\rightarrow(0,0)}
u_{ns}(t,x)
\]
such that $u_e(x)\in X^s_{e\lam}$ ($Y^e_{j\lam}$ or $Z^e_{j\lam}$
resp.). But the double limit $\lim_{(\nu,1/t)\rightarrow(0,0)}u_{ns}
(t,x)$ does not exist.

In experiment or numerical simulation of fluid motion ruled by
Navier-Stokes equations, large Reynolds number $\frac1{\nu}\gg1$ is
firstly given. Then after a long long time, we observe the fluid
motion at some time $1\ll t_1<t_2<\cdots<t_n<\cdots$. $(\nu,t_n)$ is
located in the path $\Gamma_{\mu_n}$ where $\mu_n=t_n\nu$ ($n=1,2,
\cdots$). We can see random of fluid based on the existence of path
limit of velocity along $\Gamma_{\mu_n}$ provided the initial
velocity of fluid is in one class of ${\mathbb X}_E^{2R}$,
${\mathbb X}_E$, ${\mathbb X}^s_E$, ${\mathbb Y}_E$ and
${\mathbb Z}_E$.

%%%%%%%%%%%%%%%%%%%%%%%%%%%%%%%%%%%%%%%%%%%%%%%%%

\section*{Acknowledgments}
This work is supported by National Natural Science Foundation of China--NSF,
Grant No.11971068 and No.11971077.

%%%%%%%%%%%%%%%%%%%%%%%%%%%%%%%%%%%%%%%%%%%%%%%%%%%%%%%%

\end{document}